\documentclass{amsart}
\usepackage{xspace,amssymb,mathrsfs,accents}

\newtheorem{df}{Definition}[section]
\newtheorem{prop}[df]{Proposition}
\newtheorem{theorem}[df]{Theorem}
\newtheorem{lemma}[df]{Lemma}
\newtheorem{posl}[df]{Corollary}
\newtheorem{opomba}[df]{Remark}
\newenvironment{op}%
  {\begin{opomba}\rm}{\end{opomba}}
\theoremstyle{remark}
  \newtheorem*{prooff}{Proof}

\def\theend{\hbox{}\hfill{$\surd$}}
\def\pot#1{{\mathcal #1}}

\def\mf#1{\mathfrak{#1}}
\def\ms#1{\mathscr{#1}}
\def\dom{\mathop{\rm dom}\nolimits}

\def\trcl{\mathop{\rm tr\, cl}\nolimits}
\def\cf{\mathop{\rm cf}\nolimits}
\def\Lim{\mathop{\rm Lim}\nolimits}

\def\podef{\overset{\text{def}}{\iff}}
\def\cc#1{\boldsymbol{#1}\text{\bf{-c.c.}}\xspace}

\newcommand{\forces}[2][{}]{\Vdash_{#1}\!\!\hbox{``}\,#2\,\hbox{''}}

\begin{document}

\title{Forcing $\Box_{\omega_1}$ with finite conditions}

\author{Gregor Dolinar}
\address{Biotehni\v{s}ka fakulteta\\
         Oddelek za gozdarstvo\\
         Ve\v{c}na pot 83\\
         SI - 1000 Ljubljana\\
         Slovenia}
\email{gregor.dolinar@bf.uni-lj.si}
\urladdr{}
\thanks{The research of the first author was supported in part by
        Institute of Mathematics, Physics and Mechanics, Ljubljana,
        and an ESF grant.} 

\author{Mirna D\v{z}amonja}
\address{School of Mathematics\\
         University of East Anglia\\
         Norwich, NR4 7TJ\\
         UK}
\email{m.dzamonja@uea.ac.uk}
\urladdr{http://www.mth.uea.ac.uk/~h020/}
\thanks{The research of the second author was supported by EPSRC
  grants EP/G068720 and EP/I00498X/1}

\thanks{The first author wishes to express his gratitude to University
  of East Anglia for making their facilities available during his stay
  in Norwich at the time when most of the work on this paper was done.}

\thanks{Both authors thank University of Paris VII for their
  hospitality during the final stages of the preparation of this paper.}

\thanks{This paper is a part of the first author's Ph.D. thesis in
  preparation.} 

%\date{\today}

\begin{abstract}
  We show a construction of the square principle $\Box_{\omega_1}$ by
  means of forcing with finite conditions.
\end{abstract}

\maketitle

\section{Introduction}

The square principle on a cardinal $\kappa$ states that there is a
sequence $\langle C_\alpha\rangle$ indexed by the limit ordinals in
$[\kappa,\kappa^+)$ such that each $C_\alpha$ is a club subset of 
$\alpha$ of order type $\le \kappa$ and the sequence is coherent in
the sense that if $\beta$ is a limit point of $\alpha$ then
$C_\beta=C_\alpha\cap\beta$. This principle is a feature of the
constructible universe $\mathbf L$ which was discovered by Jensen 
and used by him to show the existence of an $\omega_2$-Souslin tree in
$\mathbf L$~\cite{Jen72}. The related principle $\diamondsuit$, which
was used to construct an $\omega_1$-Souslin tree in $\mathbf L$ by Jensen,
may also be added or destroyed by forcing as wished (see~\cite{Kun83} 
for examples and discussion) and as is known of recently~(\cite{She10})
at $\kappa\ge\omega_2$ which are successors of regular cardinals, it
is simply equivalent to GCH. However, $\Box$ is connected to large
cardinals. For example, by an old proof of Solovay~\cite{SolReiKan78},
square cannot hold above a supercompact cardinal, and on smaller
cardinals, it cannot hold in the presence of forcing axioms, 
e.g.~Todor\v cevi\'c~\cite{Tod84} proved that PFA implies that for all
$\kappa\ge\omega_2$, $\Box_\kappa$ fails. Therefore $\Box$ can be seen
as a reflection principle inimical to large cardinals, and in fact by
varying the definition of square by allowing a cardinal parameter
which measures how many guesses to $C_\alpha$ we are allowed at each
$\alpha$, we obtain a hierarchy of principles of decreasing strength
which can be used to test consistency strength of various principles
(see more on this in~\cite{CumForMag01}). In the light of these facts
it is natural that the question of how to add or destroy a square
principle by forcing has been a central theme. See~\cite{CumForMag01}
for a description of some of the many known results including versions
of an older result of Magidor in which a square sequence is added by
forcing. 

The way that Magidor adds a square is to force by initial
segments along a closed bounded subsets of the domain, and to use the
existence of the ``top'' point in the domain of a forcing condition to
show that the forcing is strategically closed. Note that the principle
$\Box_\omega$ is trivially true, by taking $C_\alpha$ to be any club
of $\alpha$ of order type $\omega$, so the first non-trivial instance
of square is $\Box_{\omega_1}$. Magidor's method means that to get
$\Box_{\omega_1}$ we need to force with conditions whose domain has
size $\omega_1$. In this work we have been interested to do this
differently, using conditions whose domain is a finite set. The
interest in doing this stems from a need to understand how one can
control a one cardinal gap in forcing notions, which is a subject that
has been of interest for various combinatorial issues for a long time.
A glaring example of the need to develop this subject is the
combinatorics of $(\omega_1^{\omega_1},\le_{\rm Fin})$, which in
contrast with the vast body of knowledge about $(\omega^\omega,
\le_{\rm Fin})$ remains a mysterious object. An important development
on the subject of $(\omega_1^{\omega_1},\le_{\rm Fin})$ is Koszmider's
paper~\cite{Kos00} in which he shows that it is consistent to have an
increasing chain of length $\omega_2$ in this structure. Koszmider's
paper also gives an overview of the difficulties that there are in
forcing one gap results. 

Koszmider's method is to force with conditions where a morass is used
as a side condition. Our method is more directly connected to a
different approach, which was used to force a club on $\omega_2$ using
finite conditions. This was done in two different but similar ways by
Friedman in~\cite{Fri06} and Mitchell in~\cite{Mit09}. Both
approaches are built upon a version of adding a square on $\omega_1$
using finite conditions, as discovered by Baumgartner~\cite{Bau76} and
modified by Abraham in~\cite{AbrShe83}. The main idea in Baumgartner's
approach is that to force a club in $\omega_1$ and avoid problems at
the limit stages, one needs to specify by each condition not only what 
will go in the club, but also whole intervals that need to stay out of
it. At $\omega_2$ one can do the same, but now one needs to add side
conditions in the form of coherent systems of models in order to make
sure that cardinals are preserved, as was first done by Todor\v cevi\'c
in~\cite{Tod85}. This already is technically rather involved. What we
have done is add to this the coherent partial square sequence. Namely,
we actually force a square indexed by a club sequence -- the existence
of such a square implies the existence of an actual square sequence.
This club set is like the one added by Friedman and Mitchell. The
actual forcing notion needs to take into account the coherence of the
square sequence, and this is reflected in the complexity of the
coherence conditions between the models which form part of the forcing
conditions. An advantage of this type of approach over the morass-based
approach is that it requires less from the ground model -- for example
Friedman's forcing only needs a weakening of CH in the ground model.
We use the full CH together with $2^{\omega_1}=\omega_2$. The main
difficulties of both approaches of course are the same, and they stem
from the fact that combinatorics at $\omega_2$ is much less prone to
independence than the combinatorics at $\omega_1$, as exemplified by
the above mentioned result of Shelah on $\diamondsuit$~(\cite{She10}).
It is both in developing combinatorics and fine forcing techniques
that we can better understand the truth about $\omega_2$. 

We thank Boban Veli\v ckovi\'c for interesting discussions of
Mitchell's paper and an inspiration to consider forcing square with
finite conditions. 

\section{Preliminaries}

Most of the notation is standard. Relation $A\subset B$ means that $A$
is either a proper subset of $B$ or equal to $B$. $|X|$ is the
cardinality of set $X$. For a set of ordinals $X$, a limit point of
$X$ is an ordinal $\alpha$ such that $\alpha=\sup(Y)$ for some 
$Y\subset X$ or, equivalently, if $\alpha=\sup(X\cap\alpha)$. $\Lim(X)$
is a set of limit points of $X$. For a function $f$, $\pot D_f$
denotes the domain of $f$, and $f|_A$ denotes the restriction of $f$
to the set $A\cap\pot D_f$. If $\alpha$ and $\beta$ are ordinals then
the interval $(\alpha,\beta)$ denotes the set $\{\mu\mid\mu$ is an
ordinal, $\alpha<\mu<\beta\}=\beta\setminus(\alpha+1)$. Closed and
half open interval are defined similarly. $[A]^\kappa$ is the set of
all subsets of $A$ of cardinality $\kappa$. Set $[A]^{\le\kappa}$ is
defined analogously.

For a regular cardinal $\theta$, $H_\theta$ is the set of all sets $x$
with hereditary cardinality less than $\theta$ (i.e.~the transitive
closure of $x$ has cardinality less than $\theta$). For $\theta>\omega_2$
we consider $H_\theta$ to be a model with the standard relation $\in$
and a fixed well-ordering $\le^*$. We will primarily work with 
$H_{\omega_2}$ which we view as a model with $\in$ and 
$\le^*\upharpoonright H_{\omega_2}$. Cardinal $\theta$ is said to be
{\em large enough} if every set in  consideration is an element of
$H_\theta$.

\begin{df}
  Suppose $\kappa$ is a regular cardinal. A set $C\subset\kappa$ is
  called a {\em closed unbounded} set or a {\em club} in $\kappa$ if:

  (1) for every increasing sequence $\langle \alpha_i\mid i<\lambda
  \rangle$ of elements from $C$, for some $\lambda<\kappa$, we have
  $\bigcup_{i<\lambda}\alpha_i\in C$ {\em (closed)};

  (2) for every $\alpha<\kappa$ there exists some $\beta\in C$ such
  that $\beta>\alpha$ {\em (unbounded)}.
\end{df}

The assumption that $\kappa$ is a regular cardinal can be replaced by
a singular cardinal or even an ordinal. In that case, $\lambda$ from
clause (1) has to be below $\cf(\kappa)$. In fact, clause (1) can be
replaced by equivalent notion, that $\Lim(C)\cap\kappa\subset C$.

\begin{df}\label{DefSquare}
  Suppose $\kappa$ is a regular cardinal. {\em Square principle}
  $\Box_\kappa$ (square kappa) is a sequence 
  $\langle C_\alpha\mid\alpha$ is a limit ordinal in 
  $\kappa^+ \rangle$ such that:

  (1) $C_\alpha$ is a club in $\alpha$ for every $\alpha$;

  (2) if $\alpha\in\Lim(C_\beta)$ then $C_\alpha=C_\beta\cap\alpha$
  {\em (coherence)};
  
  (3) if $\cf(\alpha)<\kappa$ then $|C_\alpha|<\kappa$ {\em (nontriviality)}.
\end{df}

In case $\kappa=\omega_1$, the nontriviality clause simply stipulates
that if $\cf(\alpha)=\omega$ then $|C_\alpha|=\omega$.

We shall call sequence $\langle C_\alpha\mid\alpha\in\pot C\rangle$
for some set $\pot C\subset\Lim(\kappa^+)$ a {\em square-like sequence}  
if it is fulfilling all three clauses of the definition of a square
sequence. 

\section{Background on elementary submodels}

A model $M$ is an {\em elementary submodel} of a model $N$, $M\prec N$,  
if for every formula $\varphi$ with parameters $a_1,\dots,a_n\in M$,
$\varphi$ is true in $M$ if and only if it is true in $N$. If $M$ is a
countable elementary submodel of $H_\theta$ for $\theta\ge\omega_1$
then $M\cap\omega_1$ is an ordinal denoted by $\delta_M$. Also, if
$|x|\le\omega$ and $x\in M$ then $x\subset M$.

We begin by listing a few lemmas about elementary submodels which
will come in handy later on. We add proofs for completeness. When
dealing with elementary submodels, the Tarski-Vaught
test~\cite{JusWee97} comes as a very useful tool.

\begin{theorem}[Tarski-Vaught test]
  Let $M$ be a submodel of $N$. Then $M$ is elementary submodel of $N$
  if and only if for every formula $\phi(x,a_1,\dots,a_n)$ and
  $a_1\dots,a_n\in M$, if $N\models\exists x\phi(x,a_1,\dots,a_n)$
  then there exists $b\in M$ such that $N\models\phi(b,a_1,\dots,a_n)$.
\end{theorem}

\begin{lemma}\label{LemmaIntersection1}
  Suppose $N\prec H_\theta$ for some large enough $\theta$. Then
  $N\cap H_{\omega_2}\prec H_{\omega_2}$.
\end{lemma}

\begin{prooff}
  Let $a_1,\dots,a_n\in N\cap H_{\omega_2}$ and suppose that
  $H_{\omega_2}\models\psi(a_1,\dots,a_n)$ where $\psi$ is the formula
  $\exists x\phi(x,a_1,\dots,a_n)$. Then $\psi^{H_{\omega_2}}$---the
  relativization of $\psi$ to $H_{\omega_2}$---is true. Formula
  $\psi^{H_{\omega_2}}$ is equivalent to the formula $\psi^*$ obtained by
  replacing every occurence of $\exists y\in H_{\omega_2}\ \chi(y,\dots)$
  with $\exists y(\chi(y,\dots)\land|\trcl(y)|\le\omega_1)$, and
  similarly for the universal quantifier. We get $\phi^*$ from $\phi$
  in the same way. Now, $H_\theta\models\psi^*(a_1,\dots,a_n)$, or in
  other words, $H_\theta\models\exists x(\phi^*(x,a_1,\dots,a_n)
  \land|\trcl(x)|\le\omega_1)$.

  Since $\omega_1\in N$, by Tarski-Vaught test there exists some 
  $b\in N$ such that $H_\theta\models\phi^*(b,a_1,\dots,a_n)\land
  |\trcl(b)|\le\omega_1$. Hence, there exists $b\in N\cap H_{\omega_2}$
  such that $H_\theta\models\phi^{H_{\omega_2}}(b,a_1,\dots,a_n)$, and as
  a consequence, $H_{\omega_2}\models\phi(b,a_1\dots,a_n)$, which by
  Tarski-Vaught test means that $N\cap H_{\omega_2}\prec H_{\omega_2}$.
\theend
\end{prooff}

\begin{lemma}\label{LemmaIntersection2}
  Suppose $N,M\prec (H_{\omega_2},\in,\le^*)$. Then 
  $N\cap M\prec (H_{\omega_2},\in,\le^*)$.
\end{lemma}

\begin{prooff}
  Let $a_1,\dots,a_n\in N\cap M$ and suppose that $H_{\omega_2}\models
  \exists x\phi(x,a_1,\dots,a_n)$. Let $\psi(x,a_1,\dots,,a_n)$ be the
  formula $\phi(x,a_1,\dots,a_n)\land\forall y(\phi(y,a_1,\dots,a_n)\to
  x\le^*y)$. Then $H_{\omega_2}\models\exists x\psi(x,a_1,\dots,a_n)$. By  
  Tarski-Vaught test there exist $x_1\in M$ and $x_2\in N$ such that
  $H_{\omega_2}\models\psi(x_1,a_1,\dots,a_n)$ and $H_{\omega_2}\models
  \psi(x_2,a_1,\dots,a_n)$. But then $x_1=x_2=:x^*\in M\cap N$, and
  $H_{\omega_2}\models\phi(x^*,a_1,\dots,a_n)$. By Tarski-Vaught test,
  $M\cap N\prec H_{\omega_2}$.
\theend
\end{prooff}

\begin{lemma}\label{LemmaTopGap}
  If $M\prec H_\kappa$ for some $\kappa>\omega_1$, and
  $\sup(M\cap\alpha)<\alpha$ for some ordinal $\alpha\in M$, then 
  $\cf(\alpha)>\omega$.
\end{lemma}

\begin{prooff}
  If $\cf(\alpha)=\omega$ then there is a cofinal function $f:\omega
  \to\alpha$ in $M$, hence $\sup(M\cap\alpha)=\alpha$, a contradiction.
\theend
\end{prooff} 

\begin{lemma}\label{LemmaGaps}
  Let $M,N\prec H_\kappa$ for some $\kappa>\omega_1$, and suppose that
  $M\in N$. If $\alpha\not\in N$ then $\sup(M\cap\alpha)<\sup(N\cap\alpha)$.
\end{lemma}

\begin{prooff}
  If $\alpha\ge\sup(N)$ then $\sup(M\cap\alpha)=\sup(M)<\sup(N)=
  \sup(N\cap\alpha)$. Suppose now that $\alpha<\sup(N)$ and let $\beta
  :=\sup(M\cap\alpha)$ and $\beta':=\min(N\setminus\alpha)\in N$. Since 
  $M\subset N$, $\beta=\sup(M\cap\beta')$. Hence, by elementarity, 
  $\beta\in N$, and therefore $\beta<\sup(N\cap\alpha)$.
\theend
\end{prooff}

The standard reference for basic set-theoretic notions and facts
is~\cite{Jec03}. Additional source for results on elementary models in
a very concise form is~\cite{Dow88}, as well as~\cite{JusWee97}.

In our application of elementary submodels we will basically only be
interested in the ordinals that lie inside them. To simplify the
notation we will write $\ms{M}$ for a model and $M$ for its set of
ordinals $\ms{M}\cap Ord$. We will use the term ``model'' for both
$\ms M$ and $M$.
 
\section{Forcing a square}

Let $V$ be some countable transitive model of (a sufficiently large
finite fragment of) ZFC together with CH and ``$2^{\omega_1}=\omega_2$''.
Since we want to force the existence of a square sequence, the working
part of forcing notion $P$ will consist of finite partial square sequences.
We will add {\em safeguards} which will help us separate clubs from
condition $q$ and clubs from restriction $p\le q$. This will be
instrumental in the proof of properness.

It should be noted once again that we do not have to build a square
sequence on the whole $\Lim(\omega_2)$. Instead, it is enough for the
domain of the built sequence to be a club in $\omega_2$, because we
can always extend a square sequence from a club to the full 
$\Lim(\omega_2)$ (see Lemma~\ref{LemmaSquareExt}). This is the reason
why we add intervals as a part of conditions. These intervals will
serve as gaps in what will ultimately be the desired club in 
$\Lim(\omega_2)$. This way of forcing a club was introduced by 
Baumgartner in~\cite{Bau76}. 

Before we are ready to present the definition of forcing we have to
define a few auxiliary notions. For $\alpha<\omega_2$, $\cf(\alpha)=
\omega_1$, let $\ms E_\alpha$ denote some fixed countable set of clubs
in $\alpha$ of order type $\omega_1$, and $\ms E:=\langle\ms E_\alpha
\mid\alpha<\omega_2\rangle$. Define $\mf M_0:=\{\ms M\prec H_{\omega_2}
\mid\ms M$ is countable, $\ms E_\alpha\in\ms M$ for every $\alpha\in
\ms M$ with $\cf(\alpha)=\omega_1\}$. The set $\mf M_0$ will act as a
pool of possible side conditions.

For a large enough cardinal $\theta$ let $\mf M_1:=\{\ms M\prec H_\theta
\mid\ms M$ is countable, $\ms E\in\ms M\}$. Then $\mf M_1$ is a
stationary set in $[H_\theta]^\omega$. Also, if $\ms N\in\mf M_1$ and 
$\alpha\in\ms N$ then, by elementarity, $\ms E_\alpha\in\ms N$. If
$\alpha<\omega_2$ then $\ms E_\alpha\in H_{\omega_2}$, hence 
$\ms N\cap H_{\omega_2}\in\mf M_0$.

\begin{df}\label{DefCompatibleModels}
  If $\ms{M}_1,\ms{M}_2\prec H_{\omega_2}$ are countable, then we say
  that the sets $M_1$ and $M_2$ are {\em compatible} if the following
  two clauses hold as stated and with $M_1$ and $M_2$ switched:

  (a) either $M_1\cap M_2\in\ms M_1$ if $\sup(M_1\cap M_2)\in\ms M_1$,
  or $[\delta]^\omega\cap\ms M_1\subset\ms M_2$ where
  $\delta:=\sup(M_1\cap M_2)$ if $\sup(M_1\cap M_2)\not\in\ms M_1$;
  
  (b) $\{\min(M_1\setminus\lambda)\mid\lambda\in M_2,\,
  \sup(M_1\cap M_2)<\lambda<\sup(M_1)\}\cup\{\min(M_1\setminus
  \sup(M_1\cap M_2))\}$ is a finite subset of $M_1$.
\end{df}

The set in (b) is called the set of {\em $M_1$-fences for $M_2$}. This
definition of compatibility between elementary submodels (or in this
case their sets of ordinals) is due to Mitchell~\cite{Mit09}. In fact,
this version is a slight strengthening of Mitchell's compatibility
condition. The need for a slightly stronger version stems from the
fact that we have to work with sets of ordinals, namely clubs, instead
of just ordinals. Actually, workng with sets of ordinals adds a whole
new level of difficulty to the forcing construction and most of the
effort had to be invested to this end.

The following simple lemma shows that our version of compatibility
between two models is indeed stronger than Mitchell's version.

\begin{lemma}
  If $[\delta]^\omega\cap\ms M_1\subset\ms M_2$ then $M_1\cap M_2$ is
  an initial segment of $M_1$, i.e.~$M_1\cap M_2=M_1\cap\delta$.
\end{lemma}

\begin{prooff}
  Consider $\alpha\in M_1\cap\delta$, and let $A$ be some
  $\omega$-sequence of ordinals smaller than $\alpha$. Then
  $\{\alpha\}\cup A\in[\delta]^\omega\cap\ms M_1$, hence 
  $\{\alpha\}\cup A\in\ms M_2$. Therefore
  $\alpha=\max(\{\alpha\}\cup A)\in M_2$.
\theend
\end{prooff}

In lieu of the above lemma, we will say that the intersection 
$M_1\cap M_2$ is an initial segment of $M_1$ whenever clause (a) of
Definition~\ref{DefCompatibleModels} holds.

\begin{op}\label{CompatibleModels}
  If $M_1$ and $M_2$ are compatible then their structure vis-a-vis
  each other is particularly simple. Above $\sup(M_1\cap M_2)$ they
  consist of finitely many (not necessarily continuous) interchanging
  blocks, as witnessed by both fences. Below, they are either equal or
  one is a subset of the other. Namely, if there exist $\lambda_1\in
  M_1\setminus M_2$ and $\lambda_2\in M_2\setminus M_1$, $\lambda_1,
  \lambda_2<\sup(M_1\cap M_2)$, then $M_1\cap M_2$ is not an initial
  segment of either of them, hence it is an element of both $\ms M_1$
  and $\ms M_2$. Therefore, $M_1\cap M_2\in \ms M_1\cap\ms M_2$, and
  hence, $\sup(M_1\cap M_2)\in\ms M_1\cap\ms M_2$. But then 
  $\sup(M_1\cap M_2)+1\in\ms M_1\cap\ms M_2$, which is obviously a
  contradiction.
\end{op}

\begin{df}\label{DefFor}
  The forcing notion $P$ is a set of conditions of the form 
  $p:=(\pot F_p,\ \pot S_p,\ \pot O_p,\ \pot M_p)$, where
  
  (1) $\pot F_p:\Lim(\omega_2)\to\pot P(\omega_2)$, $|\pot F_p|<\omega$,
  $\pot F_p(\alpha)$ is a club $C_\alpha\subset\alpha$ of order type
  $\le\omega_1$ for all $\alpha\in\pot D_p:=\dom(\pot F_p)$, and if
  $\cf(\alpha)=\omega_1$ then $C_\alpha\in\{C\setminus\beta\mid C\in
  \ms E_\alpha,\,\beta\in\pot D_p\cap\alpha\}$; 

  (2) $\pot S_p\subset\pot D_p$ and $\alpha\in\pot S_p$ for every
  $\alpha\in\pot D_p$ with $\cf(\alpha)=\omega_1$;

  (3) $\pot M_p$ is a finite set such that if $M\in\pot M_p$ then
  there exists a countable elementary submodel $\ms M\in\mf M_0$
  with $M=\ms M\cap Ord$, and additionally, $\sup(M)\in\pot S_p$ for
  every $M\in\pot M_p$;

  (4) for every $\alpha,\beta\in\pot D_p$, $\alpha\ne\beta$, if
  $\mu\in\Lim(C_\alpha)\cap\Lim(C_\beta)$ then $C_\alpha\cap\mu=C_\beta\cap\mu$;

  (5) if $\alpha\in\pot D_p$ and $\sigma\in\pot S_p$, $\sigma<\alpha$,
  then $C_\alpha\cap\sigma$ is a finite set; 

  (6) for all $\alpha\in\pot D_p$ and $M\in\pot M_p$:
  
  (a) if $\alpha\in M$ then $C_\alpha\in\ms M$,

  (b) if $\alpha\not\in M$ is such that $\alpha<\sup(M)$, or if $\alpha\in M$ 
  is such that $\sup(M\cap\alpha)<\alpha$, then $\min(M\setminus\alpha)
  \in\pot S_p$ and $\sup(M\cap\alpha)\in\pot D_p$,

  (c) if $\alpha\not\in M$ and $\sup(M\cap\alpha)<\alpha<\sup(M)$ then
  $C_\alpha=C_\beta\cap\alpha$ if there is some $\beta\in\pot D_p$,
  $\beta>\alpha$, such that $\alpha\in\Lim(C_\beta)$, otherwise 
  $C_\alpha\cap\sup(M\cap\alpha)$ is a finite set,

  (d) if $\alpha\not\in M$ and $\sup(M\cap\alpha)=\alpha$ then 
  $C_\alpha=C_\beta\cap\alpha$ if there is some $\beta\in\pot D_p$,
  $\beta>\alpha$, such that $\alpha\in\Lim(C_\beta)$, otherwise
  $C_\alpha$ is some cofinal sequence in $\alpha$ of length $\omega$;

  (7) $\pot O_p$ is a finite set of half open nonempty intervals 
  $(\beta',\beta]\subset\omega_2$ such that $\pot D_p
  \cap\bigcup\pot O_p=\emptyset$; 

  (8) if $(\beta',\beta]\in\pot O_p$ and $M\in\pot M_p$
  then either $(\beta',\beta]\in\ms M$ or 
  $(\beta',\beta]\cap\ms M=\emptyset$; 

  (9) if $M_1,M_2\in\pot M_p$ then they are compatible, and the
  $M_1$-fence for $M_2$ and the $M_2$-fence for $M_1$ are both subsets
  of $\pot S_p$.\\[5pt]  
\noindent
  For $p,q\in P$ define $p\le q\podef\pot F_p\subset\pot F_q,\,\pot S_p
  \subset \pot S_q,\,\pot O_p\subset\pot O_q,\,\pot M_p\subset\pot M_q$.
\end{df}

Clause (6b) tells us that a gap in a model $M$ has to be closed from above
by a safeguard if there is something (i.e.~an ordinal $\alpha\in\pot D_p$)
inside that gap. This safeguard is an echo of $\alpha$ resonating in
$M$, warning everybody in $M$ to stay away from that gap. Fences from
clause (9) serve exactly the same purpose.

Notice that in clause (8), the interval $(\beta',\beta]$ is an element
of the model $\ms M$ if and only if both $\beta'$ and $\beta$ are in $M$.

\begin{lemma}
  $(P,\le)$ is a separative forcing notion.
\end{lemma}

\begin{prooff}
  Transitivity is trivial. The minimal element is
  $(\emptyset,\emptyset,\emptyset,\emptyset)$. For separativeness
  consider an arbitrary condition $p\in P$. We will find two incompatible
  extensions. Let $\alpha:=\sup(\pot D_p\cup\bigcup\pot O_p\cup
  \bigcup\pot M_p)$, and $\beta:=\alpha+\omega<\omega_2$. Define
  $C_\beta:=[\alpha,\beta)$ and $C'_\beta:=(\alpha,\beta)$. It is easy
  to check that $q:=(\pot F_p\cup\{(\beta,C_\beta)\},\,\pot S_p,\,
  \pot O_p,\,\pot M_p)$ and $q':=(\pot F_p\cup\{(\beta,C'_\beta)\},\,
  \pot S_p,\,\pot O_p,\,\pot M_p)$ are both conditions extending $p$,
  and that they are incompatible. Notice, that since $\cf(\beta)=\omega$,
  $C_\beta$ and $C'_\beta$ need not be in $\ms E_\beta$.
\theend
\end{prooff}

We first prove several lemmas that show us a little bit more about the
structure of conditions in $P$, and will be helpful in further proofs.
Most notably, they will shed some lights on the correspondence between
models and clubs, and thus clarify clause (6).

\begin{lemma}\label{LemmaLimPoint1}
  Let $p\in P$, and suppose that $\alpha,\gamma\in\pot D_p$ and
  $M\in\pot M_p$ are such that $\alpha<\sup(M)$, $\alpha\not\in M$, and
  $\alpha\in\Lim(C_\gamma)$. Then $\gamma\le\min(M\setminus\alpha)$.
\end{lemma}

\begin{prooff}
  By (6b), $\sigma:=\min(M\setminus\alpha)\in\pot S_p$, hence if
  $\gamma>\sigma$ then, by (5), $C_\gamma$ has no limit points below
  $\sigma$.
\theend
\end{prooff}

Notice that if $\alpha\in\Lim(C_\gamma)$ then $\cf(\alpha)=\omega$,
otherwise $C_\gamma$ would have order type larger than $\omega_1$.

\begin{lemma}\label{LemmaLimPoint2}
  If $p,\alpha,\gamma$ and $M$ are as in the previous lemma and $\max
  \{\gamma\in\pot D_p\mid\alpha\in\Lim(C_\gamma)\}<\min(M\setminus\alpha)$,
  or if $\alpha>\sup(M)$, then $C_\alpha\cap\sup(M\cap\alpha)$ is finite,
  and therefore $\sup(M\cap\alpha)<\alpha$.
\end{lemma}

\begin{prooff}
  If $\alpha>\sup(M)$ then the conclusion follows from clauses (3) and
  (5). Suppose now that $\alpha<\sup(M)$ and $\alpha\in\Lim(C_\gamma)$.  
  Then $\gamma\not\in M$ and $\sup(M\cap\gamma)>\gamma$. If $\gamma$
  is not a limit point of any $C_{\gamma'}$ then, by (6c), 
  $C_\gamma\cap\sup(M\cap\gamma)$ is finite. Since $C_\alpha\subset C_\gamma$
  and $\sup(M\cap\alpha)=\sup(M\cap\gamma)$, $C_\alpha\cap\sup(M\cap\alpha)$
  is also finite. If $\gamma\in\Lim(C_{\gamma'})$ then 
  $\alpha\in\Lim(C_{\gamma'})$ and we can repeat the above argument. As  
  $|\pot D_p|<\omega$, we only have to repeat it finitely many times,
  and in the end we can conclude that $C_\alpha\cap\sup(M\cap\alpha)$
  is finite.
\theend
\end{prooff}

\begin{lemma}\label{LemmaLimPoint3}
  If $M\in\pot M_p$ is countable then $C_{\sup(M)}$ is an $\omega$-sequence.
\end{lemma}

\begin{prooff}
  By clauses (3) and (5), $\sup(M)$ cannot be a limit point of any
  $C_\gamma$ for $\gamma\in\pot D_p$. So $C_{\sup(M)}$ is an
  $\omega$-sequence by clause (6d).
\theend
\end{prooff}

\begin{lemma}\label{LemmaQ}
  Let $\ms N'\in\mf M_1$ be a countable elementary submodel of $H_\theta$.
  If $p$ is a condition in $P\cap\ms N'$ then there exists an extension 
  $q\ge p$ such that $\ms N'\cap H_{\omega_2}\in\pot M_q$.
\end{lemma}

\begin{prooff}
  Let $p$ be of the form $(\pot F_p,\,\pot S_p,\,\pot O_p,\,\pot M_p)$
  and let $\ms N:=\ms N'\cap H_{\omega_2}\in\mf M_0$. By
  Lemma~\ref{LemmaIntersection1}, $\ms N\prec H_{\omega_2}$. Note also
  that $p\in\ms N$.

  For every $\alpha\not\in N$ such that $\alpha=\sup(N\cap\gamma)$ for
  some $\gamma\in\pot D_p$, let $C_\alpha$ be a club according to
  clause (6d). In the case of $\alpha\in\Lim(C_\beta)$ for some
  $\beta\in\pot D_p$ this choice is well-defined because by clause (4)
  it does not depend on $\beta$. It is worth mentioning that
  $\cf(\alpha)=\omega$, hence $C_\alpha$ need not be in $\ms E_\alpha$. 
  Notice, that by Lemma~\ref{LemmaTopGap}, $\cf(\gamma)=\omega_1$.
  Therefore, $\gamma$ is already in $\pot S_p$ and does not have to be
  added for clause (6b) to be satisfied. We will also have to add
  $\sup(N)$ to the set of safeguards. For the corresponding club
  $C_{\sup(N)}$ we pick any cofinal $\omega$-sequence. Again,
  $\cf(\sup(M))=\omega$, therefore $C_{\sup(M)}$ des not have to be in
  $\ms E_{\sup(M)}$.

  Define $q:=(\pot F_p\cup\{(\alpha,C_\alpha)\mid\alpha\not\in N,\,
  \alpha=\sup(N\cap\gamma)$ for some $\gamma\in\pot D_p\}\cup
  \{(\sup(N),C_{\sup(N)})\},\,\pot S_p\cup\{\sup(N)\},\,\pot O_p,\,
  \pot M_p\cup\{N\})$. Most of the clauses of definition~\ref{DefFor}
  are trivially true, including clause (6), which is due to the fact
  that we used clause (6) to construct additional clubs. These new clubs
  conform to clauses (1) and (4) as well, since clause (6d) was added
  specifically for this purpose. Notice, that for $M\in\pot M_p$ the
  $M$-fence for $N$ is the empty set, while the $N$-fence for $M$ is
  $\{\sup(M\cap N)\}=\{\sup(M)\}$ which is a subset of $\pot S_p\subset
  \pot S_q$ by clause (3).  As for clause (7), suppose that some newly
  added $\alpha<\sup(N)$ falls into some interval $(\beta',\beta]$.
  Then its corresponding $\gamma\in\pot D_p$ was already in this
  interval, since $\{\beta',\beta\}\subset N$. But that is in a
  contradiction with clause (7) in $p$.

  Hence $q$ is the desired condition extending $p$. 
\theend
\end{prooff}

\begin{lemma}\label{LemmaR}
  Let $\ms N\in\mf M_1$ be a countable elementary submodel of
  $H_{\omega_2}$, and suppose that $r\in P\setminus\ms N$ is such that
  $N\in\pot M_r$. Then $r_{\ms N}:=(\pot F_r\cap\ms N,\,
  \pot S_r\cap\ms N,\,\pot O_r\cap\ms N,\,
  (\pot M_r\cap\ms N)\cup\{M\cap N\mid M\in\pot M_r,\,M\not\in\ms N,\,
  M\cap N\in\ms N\})$ is a condition in $P\cap\ms N$.
\end{lemma}

\begin{prooff}
  First note that by clause (6a), $\pot D_{r_{\ms N}}=\pot D_r\cap\ms N$.
  Also, by Lemma~\ref{LemmaIntersection2}, $\ms M\cap\ms N\prec
  H_{\omega_2}$. Additionally, $\ms M\cap\ms N\in\mf M_0$, hence
  $M\cap N$ can be added to $\pot M_{r_{\ms N}}$. Notice that if
  $M\cap N\in\ms N$ then $\sup(M\cap N)\in\pot S_r\cap\ms N$ because
  it is in the $N$-fence for $M$, hence clause (3) is satisfied.

  It is obvious that $r_{\ms N}\in\ms N$. Slightly less trivial thing
  to prove is that $r_{\ms N}\in P$. Compatibility between the
  elements of $\pot F_r$, $\pot S_r$, $\pot O_r$ and $\pot M_r$ is
  inherited from $r$, as is the compatibility between two models from
  $\pot M_r$. The same can be said for the compatibility between
  $\pot F_{r_{\ms N}}$, $\pot O_{r_{\ms N}}$ and a model of the form
  $M\cap N\in\pot M_{r_{\ms N}}$, however with a closer inspection of
  clause (6). 

  For clause (6b) consider $\alpha\in\pot D_{r_{\ms N}}$ and $M\cap N\in
  \pot M_{r_{\ms N}}$ such that $\alpha\not\in M\cap N$. That means
  that $\alpha\not\in M$. Since $M\cap N\in\ms N$, $M\cap N$ is an
  initial segment of $M$. If $\alpha<\sup(M\cap N)$ then $\min((M\cap N)
  \setminus\alpha)=\min(M\setminus\alpha)\in\pot S_r\cap\ms N$ by
  clause (6b) in $r$, hence $\min((M\cap N)\setminus\alpha)\in
  \pot S_{r_{\ms N}}$. By the same argument, $\sup((M\cap N)\cap\alpha)\in
  \pot D_{r_{\ms N}}$. Similarly, if $\alpha\in M\cap N$ such that 
  $\sup((M\cap N)\cap\alpha)<\alpha$ then by (6b) in $r$, $\alpha=
  \min((M\cap N)\setminus\alpha)\in\pot S_{r_{\ms N}}$
  and $\sup((M\cap N)\cap\alpha)\in\pot D_{r_{\ms N}}$. 

  For clause (6c) assume that $\alpha\not\in M\cap N$ is such that
  $\sup((M\cap N)\cap\alpha)<\alpha<\sup(M\cap N)$. The only potential
  problem is if $\alpha\in\Lim(C_\beta)$ for some $\beta\in\pot D_r$
  while $\alpha\not\in\bigcup_{\gamma\in\pot D_{r_{\ms N}}}\Lim(C_\gamma)$. In
  such a case, by Lemma~\ref{LemmaLimPoint1}, $\beta\le
  \min(M\setminus\alpha)$. In fact, $\beta<\min(M\setminus\alpha)$
  because $\min(M\setminus\alpha)\in N$. Then by
  Lemma~\ref{LemmaLimPoint2}, $C_\alpha\cap\sup(M\cap\alpha)$ is finite.

  For clause (6d) assume that $\alpha$ is a supremum of a block of
  $M\cap N$. Let us ask ourselves a question. Is it possible that
  $\alpha$ is a limit point of some $C_\beta$ in $r$ but not in
  $r_{\ms N}$? It is only possible if $\alpha=\sup(M\cap N)$. Otherwise, 
  as above, $\min(M\setminus\alpha)\in\pot S_r$ which makes it the only 
  candidate for $\beta$. However, this $\beta$ and $C_\beta$ remain in
  $r_{\ms N}$ and do not answer our question affirmatively. On the
  other hand $\sup(M\cap N)\in\pot S_r$ by (9) in $r$ and hence by (5)
  cannot be a limit point of any $C_\beta$, neither in $r$ nor in
  $r_{\ms N}$.

  The thing that merits the closest attention is the compatibility
  between two models of the form $M\cap N\in\pot M_{r_{\ms N}}$. Consider
  $M_i':=M_i\cap N$, $i=1,2$, such that $M_i\in\pot M_r\setminus\ms N$
  and $M_i'\in\ms N$. Let $x_1$ be the $M_1$-fence for $M_2$. Then
  $x_1\cap N=x_1\cap\sup(M_1\cap N)$ is the $M_1'$-fence for $M_2'$.
  Here we use the fact that $M_1\cap N$ is an initial segment of
  $M_1$. We get the $M_2'$-fence for $M_1'$ in a similar way.

  Suppose that $M_1\cap M_2\in\ms M_2$. Let $\lambda:=\sup(M_1'\cap
  M_2')\in N$. We have to consider three cases.

  {\it Case 1:}  $\lambda\in M_2$. Then, by elementarity,
  $M_1'\cap M_2'=M_1\cap M_2\cap\lambda\in\ms M_2\cap\ms N=\ms M_2'$.

  {\it Case 2:} $\lambda\not\in M_2$ and $\lambda\in M_1$. Obviously
  $\lambda\le\sup(M_1\cap M_2)$. If $\lambda<\sup(M_1\cap M_2)$ then
  $\lambda\in M_2$, because $M_1\cap M_2$ is an initial segment of
  $M_1$. Hence we get a contradiction. If $\lambda=\sup(M_1\cap M_2)$
  then $\lambda\in M_2$, because $M_1\cap M_2\in\ms M_2$, and again
  we get a contradiction.

  {\it Case 3:} $\lambda\not\in M_2$ and $\lambda\not\in M_1$. If
  $\lambda=\sup(M_1\cap M_2)$ then, as above, $\lambda\in M_2$. If
  $\lambda<\sup(M_1\cap M_2)$ then $\lambda=\sup((M_1\cap M_2)\cap\min
  (M_2\setminus\lambda))\in M_2$ by elementarity. Both possibilities
  lead to a contradiction. 

  Suppose now that $[\mu]^\omega\cap\ms M_2\subset\ms M_1$ where
  $\mu:=\sup(M_1\cap M_2)$. If $A\in[\lambda]^\omega$ is in $\ms M_2'$
  then $A\in\ms M_2\cap\ms N$ and $A\in[\mu]^\omega$, hence
  $A\in\ms M_1\cap\ms N=\ms M_1'$. So $[\lambda]^\omega\cap\ms M_2'
  \subset\ms M_1'$.

  Basically, $M_1'\cap M_2'$ is in the same correspondence with both
  $\ms M_1'$ and $\ms M_2'$ as $M_1\cap M_2$ is with $\ms M_1$ and
  $\ms M_2$ respectively.
\theend
\end{prooff}

We are now ready to prove the most important facet of forcing $P$,
namely the fact that it preserves $\omega_1$. We do that by proving
that $P$ is proper. There are several equivalent definitions of
properness. We shall use the following one.

\begin{df}\label{DefProperness}
  Let $P$ be a forcing notion, $\theta$ a large enough cardinal, and
  let $\mf N$ be some set of countable elementary submodels of $H_\theta$,
  stationary in $[H_\theta]^\omega$.

  (1) Condition $q\in P$ is {\em $\ms N$-generic} if for every
  extension $r\ge q$, $r\in P$, and every dense set $\ms D\subset P$,
  $\ms D\in\ms N$, there exists some condition $s\in\ms D\cap\ms N$
  which is compatible with $r$.

  (2) $P$ is {\em proper} if for every $\ms N\in\mf N$ such that
  $P\in\ms N$, every condition $p\in P\cap\ms N$ has an $\ms N$-generic
  extension.
\end{df}

\begin{prop}
  Forcing $P$ is proper.
\end{prop}

\begin{prooff}
  Let $\theta$ be a large enough cardinal. Fix a countable elementary
  submodel $\ms N'\prec H_\theta$, $\ms N'\in\mf M_1$, such that
  $P\in\ms N'$, and consider an arbitrary $p=(\pot F_p$, $\pot S_p$,
  $\pot O_p$, $\pot M_p)\in P\cap\ms N'$. Define $\ms N:=\ms N'\cap
  H_{\omega_2}\in\mf M_0$ and let $q$ be an extension of $p$ given by
  Lemma~\ref{LemmaQ}. We will prove that $q$ is an 
  $\ms N'$-generic extension of $p$.

  Suppose $r\in P$ is an arbitrary extension of $q$. Let $r_{\ms N}$
  be the condition given by Lemma~\ref{LemmaR}. Proceed by fixing a
  dense open subset $\ms D\subset P$, $\ms D\in\ms N'$, and extend
  $r_{\ms N}$ to $s\in\ms D\cap\ms N'$. Clearly $s\in H_{\omega_2}$, hence 
  $s\in\ms N$. We shall prove that $r$ and $s$ are compatible by
  proving clause by clause of Definition~\ref{DefFor} that (a certain
  extension of) $t:=(\pot F_r\cup\pot F_s,\,\pot S_r\cup\pot S_s,\,
  \pot O_r\cup\pot O_s,\,\pot M_r\cup\pot M_s)$ is a condition.

  {\it Clauses (1), (2) and (3)} are obviously true.

  {\it Clause (4):} take arbitrary $\alpha,\beta\in\pot D_t$. We can
  assume WLOG that $\alpha\in\pot D_r\setminus\pot D_s$ and $\beta\in
  \pot D_s\setminus\pot D_r$. If $\beta>\alpha$ then there are two
  possibilities. If $\beta=\min(N\setminus\alpha)$ then $\beta\in\pot
  S_r\subset\pot D_r$ by (6b) in $r$, which we assumed was not the
  case. If $\beta>\min(N\setminus\alpha)$ then $C_\beta\cap C_\alpha
  \subset C_\beta\cap\min(N\setminus\alpha)$ which is finite by (5) in
  $s$, because $\min(N\setminus\alpha)\in\pot S_r\cap\ms N\subset \pot
  S_s$ by (6b) in $r$. Hence, $\Lim(C_\beta)\cap\Lim(C_\alpha)=\emptyset$.

  If $\beta<\sup(N)\le\alpha$ then $C_\alpha\cap\beta$ is finite
  either by Lemma~\ref{LemmaLimPoint3} or clause (5) in $r$, hence
  $\Lim(C_\alpha)\cap\Lim(C_\beta)=\emptyset$. If $\beta<\alpha<\sup(N)$, 
  then by clause (6c) or (6d) in $r$, either $C_\alpha\cap\delta$ is
  finite for every $\delta<\sup(N\cap\alpha)$ including $\beta$, or
  $C_\alpha=C_{\gamma'}\cap\alpha$ for some $\gamma'\in\pot D_r$. In
  the latter case let $\gamma$ be the largest such $\gamma'$. By
  Lemma~\ref{LemmaLimPoint1}, $\gamma\le\min(N\setminus\alpha)$. If
  $\gamma<\min(N\setminus\alpha)$ then $C_\alpha\cap\beta$ is finite
  by Lemma~\ref{LemmaLimPoint2}, because $\beta<\sup(N\cap\alpha)$. If
  $\gamma=\min(N\setminus\alpha)\in\pot D_s$ and if $C_\alpha$ and
  $C_\beta$ have a common limit point $\mu$, then $\mu\in\Lim(C_\gamma)$,
  hence $C_\beta\cap\mu=C_\gamma\cap\mu=C_\alpha\cap\mu$. The first
  equality follows from (4) in $s$ and the second follows from (4) in
  $r$.

  If $\sup(N\cap\alpha)=\alpha$ then the argument is similar. If
  $C_\alpha$ is an $\omega$-sequence of ordinals then (4) is trivially
  true. If $C_\alpha=C_\gamma\cap\alpha$ for some $\gamma\in\pot D_r$
  and if $\alpha<\sup(N)$, then $\gamma=\min(N\setminus\alpha)$.
  Notice that if $\gamma$ were below $\min(N\setminus\alpha)$ then
  $C_\alpha$ would be finite by (6c) in $r$. Now we get (4) in $t$
  just as above. If $\alpha=\sup(N)$ then $C_\alpha$ were again finite
  by (6c) or even (5) in $r$. 

  {\it Clause (5):} first consider arbitrary $\alpha\in\pot D_r\setminus
  \pot D_s$ and $\sigma\in\pot S_s\setminus\pot S_r$, $\sigma<\alpha$.
  By (6c) or (6d) in $r$, either $C_\alpha=C_\gamma\cap\alpha$ for some
  $\gamma\in\pot D_r$, or $C_\alpha\cap\delta$ is a finite set for
  every $\delta<\sup(N\cap\alpha)$. The second case is trivial,
  because $\sigma\in N$, hence $\sigma<\sup(N\cap\alpha)$. In the
  first case, if $\alpha<\sup(N)$ then we can assume that $\gamma=
  \min(N\setminus\alpha)\in\pot D_s$. Hence, $C_\alpha\cap\sigma=
  C_\gamma\cap\sigma$ which is a finite set by (5) in $s$. If $\alpha>
  \sup(N)$ then $C_\alpha\cap\sup(N)$ is finite regardless of
  $\gamma$, and if $\alpha=\sup(N)$ then $C_\alpha$ is an
  $\omega$-sequence.

  If $\alpha\in\pot D_s\setminus\pot D_r$ and $\sigma\in\pot S_r
  \setminus\pot S_s$, $\sigma<\alpha$, then $\min(N\setminus\sigma)\in
  \pot S_r\cap\ms N\subset\pot S_s$. Also, $\alpha\ge\min(N\setminus
  \sigma)$, but $\alpha\ne\min(N\setminus\sigma)$, otherwise $\alpha\in
  \pot S_r\subset\pot D_r$. Hence $\alpha>\min(N\setminus\sigma)
  \in\pot S_s$ and therefore $C_\alpha\cap\sigma\subset C_\alpha\cap
  \min(N\setminus\sigma)$ which is a finite set by (5) in $s$.

  {\it Clause (6):} first consider arbitrary $\alpha\in\pot D_r\setminus
  \pot D_s$ and a $M\in\pot M_s\setminus\pot M_r$. Then $\alpha\not
  \in N\supset M$ and $\sup(M\cap\alpha)<\alpha$ by Lemma~\ref{LemmaGaps}.
  We argue just like for clause (5). Either $C_\alpha=C_\gamma\cap
  \alpha$ for some $\gamma\in\pot D_r\subset\pot D_t$, and we get (6c)
  in $t$, or $C_\alpha\cap\delta$ is a finite set for every $\delta<
  \sup(N\cap\alpha)$. Since $\sup(M\cap\alpha)<\sup(N\cap\alpha)$, we
  also get (6c).

  Additionally, if $\alpha<\sup(M)$ then let $\sigma:=\min(N\setminus
  \alpha)\in\pot D_s$. Since $M\subset N$, we have $\sup(M\cap\sigma)<
  \sigma$, hence, by (6b) in $s$, $\sup(M\cap\alpha)=\sup(M\cap\sigma)
  \in\pot D_s\subset\pot D_t$, and $\min(M\setminus\alpha)=
  \min(M\setminus\sigma)\in\pot S_s\subset\pot S_t$. 

  Now consider arbitrary $\alpha\in\pot D_s\setminus\pot D_r$ and
  $M\in\pot M_r\setminus\pot M_s$. Assume first that $\alpha
  \not\in M$. If $\sup(M\cap N)<\alpha<\sup(M)$ then $\alpha':=
  \min(M\setminus\alpha)\in\pot S_r\subset\pot S_t$ by (9) in $r$
  because there is a block of $N$ inside a gap of $M$. Also, by (6b)
  in $r$, $\sup(M\cap\alpha)=\sup(M\cap\alpha')\in\pot D_r\subset
  \pot D_t$. Additionally (and also if $\alpha>\sup(M)$), if there is
  no $\delta\in\pot D_t$ such that $\alpha\in\Lim(C_\delta)$ then
  $C_\alpha\cap\sup(M\cap\alpha)$ is finite, because either there is
  some $\gamma\in M$, $\sup(M\cap N)<\gamma<\alpha$, in which case
  $\sigma:=\min(N\setminus\gamma)\in\pot S_r\cap N\subset\pot S_s$ by
  (9) in $r$, or else $\sigma:=\min(N\setminus\sup(M\cap N))\in\pot S_r 
  \cap N$, because $\sigma$ is in the $N$-fence for $M$. In both cases
  $C_\alpha\cap\sup(M\cap\alpha)\subset C_\alpha\cap\sigma$ which is
  finite by (5) in $s$. Hence, we got (6c) which we had to get, because
  $\sup(M\cap\alpha)<\alpha$. Notice, that if $\sup(M\cap\alpha)=
  \alpha$ then $\alpha$ is a minimum of a block in $N$ and consequently
  $\alpha\in\pot S_r\subset\pot D_r$ by (9) in $r$, which is also the
  case if $\alpha=\sup(M\cap N)$.

  However, if $\alpha<\sup(M\cap N)$ then $M\cap N\in\ms N$. Therefore,
  by compatibility between $\alpha$ and $M\cap N$ in $s$, $\min(M
  \setminus\alpha)=\min((M\cap N)\setminus\alpha)\in\pot S_s\subset
  \pot S_t$ and $\sup(M\cap\alpha)=\sup((M\cap N)\cap\alpha)\in\pot D_s
  \subset\pot D_t$. Also, if there is no $\delta\in\pot D_t$ such that
  $C_\alpha\in\Lim(C_\delta)$ then $C_\alpha\cap\sup(M\cap\alpha)=C_\alpha 
  \cap\sup((M\cap N)\cap\alpha)$ which is either an $\omega$-sequence if
  $\sup(M\cap\alpha)=\alpha$, or finite if  $\sup(M\cap\alpha)<\alpha$.
  If $\alpha\in\Lim(C_\delta)$ for some $\delta\in\pot D_s$ then
  $C_\alpha=C_\delta\cap\alpha$ by (6) in $s$. If such $\delta$ is in
  $\pot D_r\setminus\pot D_s$ then, by (4) which was already proved
  for $t$, $C_\alpha=C_\delta\cap\alpha$.

  If $\alpha\in M$ then we have to consider two possibilities. The
  first case is if $M\cap N\in\ms N$. Then, by (6a) in $s$, $C_\alpha
  \in\ms M$. Also, if $\sup(M\cap\alpha)<\alpha$ then $\sup(M\cap\alpha)=
  \sup((M\cap N)\cap\alpha)\in\pot D_s\subset\pot D_t$, and $\min(M
  \setminus\alpha)=\min((M\cap N)\setminus\alpha)\in\pot S_s\subset
  \pot S_t$, both by (6b) in $s$.

  If, on the other hand, $M\cap N$ is an initial segment of
  $N$, and if $\cf(\alpha)=\omega$ then, since $C_\alpha\in
  [\sup(M\cap N)]^\omega\cap\ms N$, by
  Definition~\ref{DefCompatibleModels}, $C_\alpha\in\ms M$. If
  $\cf(\alpha)=\omega_1$ then $C_\alpha=C\setminus\beta$ for some
  $C\in\ms E_\alpha\in\ms M$ and $\beta\in\pot D_s\cap\alpha\subset
  \ms N$. Since $\beta<\alpha$ and $M\cap N$ is an initial segment of
  $N$, we have $\beta\in\ms M$, hence $C_\alpha\in\ms M$. Note that if
  $\pot D_s\cap\alpha=\emptyset$ then $C_\alpha$ is in $\ms E_\alpha$,
  and by extension in $\ms M$, by itself.

  We also have to consider the possibility that $\sup(M\cap\alpha)<
  \alpha$. Then $\sup(M\cap\alpha)\ge\sup(N\cap\alpha)$. Here we are
  in trouble because $\alpha$ knows nothing about $M$, so there is no
  reason to believe that $\sup(M\cap\alpha)$ is either in $\pot D_r$
  or in $\pot D_s$. In fact, it is certainly not in $\pot D_r$ because
  then $\alpha$ would also be in $\pot D_r$ by (6b) in $r$ applied on
  $M$ and $\sup(M\cap\alpha)$. This is the reason why we have to
  enlarge $\pot D_t$ with every possible $\beta:=\sup(M\cap\alpha)$
  for all pairs $\alpha\in\pot D_s\setminus\pot D_r$ and countable
  $M\in\pot M_r\setminus\pot M_s$ such that $M\cap N\not\in\ms N$,
  $\alpha<\sup(M\cap N)$ and $\sup(M\cap\alpha)<\alpha$. Notice, that
  this way we also prove the consistency between $\alpha$ and $N$.
  $C_{\beta}$'s have to be chosen according to (6d) in $t$. The
  precise definition will follow. Meanwhile, it is worth mentioning
  that $\beta\not\in\pot D_t$, hence $C_{\beta}$ is a new club, and we
  are not in danger of overwriting an existing club with a new and
  different one. If $\beta$ were present in $t$ it could only come
  from $r$, as $\beta\not\in N$. But then $\alpha$ itself would be
  present by (6b) in $r$. Also, $\cf(\beta)=\omega$, hence we need not
  worry about $\ms E_\beta$.

  We have to check that every new $\beta$ and its club are compatible
  with everything in $t$. Club $C_{\beta}$ is compatible with all the
  other clubs as well as with $\pot S_t$ because of its construction.
  Point $\beta$ is compatible with $\pot O_t$ because $\alpha$ is --
  at least it will be, once we prove clause (7) in $t$. The only
  nontrivial thing to check is the compatibility between $\beta$ and
  $M'\in\pot M_r\setminus\pot M_s$. If $M'\in\pot M_s$ then the
  compatibility follows from compatibility between $\alpha$ and $M'$
  in $s$. Similarly, we only have to consider $M'$ such that
  $M'\cap N\not\in\ms N$, otherwise we simply use either the compatibility
  between $M'\cap N$ and $\alpha$ in $s$ if $\alpha<\sup(M'\cap N)$,
  or the compatibility between $M'$ and $N$ in $r$ if $\alpha\ge\sup
  (M'\cap N)$.

  Suppose that $M'$ is as described. If $\sup(M\cap M')<\beta<\sup(M')$
  then the block of $M$ with supremum $\beta$ lies in a gap of $M'$.
  But then, by (9) in $r$, $\min(M'\setminus\beta)\in\pot S_r$, and
  consequently $\sup(M'\cap\beta)\in\pot D_r$ by (6b) in $r$. In fact,
  in this case $\beta$ was already present in $r$ and did not have to
  be added. If $\beta=\sup(M\cap M')$ then $\beta\in\pot D_r$, hence
  $\alpha=\min(M\setminus\beta)\in\pot S_r$.

  Suppose now that $\beta<\sup(M\cap M')$. If $M\cap M'\in\ms M'$
  then $\beta\in M'$ and it is not even a minimum of a block. If
  $C_{\beta}=C_\alpha\cap\beta$ then $C_{\beta}\in\ms M'$, because
  $\alpha\in M'$. However, if $C_{\beta}$ is an $\omega$-sequence then 
  it has to be chosen so that it is an element of every such $\ms M'$.
  Therefore, the precise definition of $C_\beta$ is this: if
  $\beta\in C_\alpha$ then $C_\beta:=C_\alpha\cap\beta$, otherwise
  $C_\beta$ is some $\omega$-sequence cofinal in $\beta$ such that
  $C_\beta\in\bigcap\{M'\in\pot M_r\setminus\pot M_s\mid\beta\in M'\}
  \prec H_{\omega_2}$. In the second case, we can find such
  $\omega$-sequence by elementarity, since a finite intersection of
  elementary submodels of $H_{\omega_2}$ is itself an elementary
  submodel of $H_{\omega_2}$.

  If $M\cap M'\in\ms M$ and $\alpha\le\sup(M'\cap N)$ then, because
  $M'\cap N$ is an initial segment of $N$, $\alpha$ is also a minimum
  of a block in $M'$ and hence $\sup(M'\cap\beta)=\sup(M'\cap\alpha)$, 
  which had to be added to $\pot D_t$ just like $\beta$. Additionally,
  $\beta$ is inside a gap of $M'$ and $C_\beta$ is compatible with
  $M'$ because it was constructed in the spirit of (6d). Now, if
  $\alpha>\sup(M'\cap N)$ then a block of $N$ with $\alpha$ in it lies
  in a gap of $M'$. But then $\sup(M'\cap\beta)=\sup(M'\cap\alpha)\in
  \pot D_r\subset\pot D_t$ and $\min(M'\setminus\beta)=\min(M'
  \setminus\alpha)\in\pot S_r\subset\pot S_t$. Again, $\beta$ is
  inside a gap of $M'$ and $C_\beta$ is compatible with $M'$ because
  of its construction.

  If $M\cap M'$ is an initial segment of both $M$ and $M'$ then they
  have the same gaps below $\sup(M\cap M'$) and $\beta$ closes the
  gap in $M'$ just as it does it in $M$. Also, $C_{\beta}$ is then
  compatible with $M'$ just as it is with $M$.

  {\it Clause (7):} first consider arbitrary $(\beta',\beta]\in\pot O_r
  \setminus\pot O_s$ and $\alpha\in\pot D_s\setminus
  \pot D_r$. Since $(\beta',\beta]\cap N=\emptyset$ by
  (8) in $r$, and $\alpha\in N$, we can conclude that $\alpha\not\in
  (\beta',\beta]$. Now, if $\alpha\in\pot D_r\setminus\pot D_s$
  and $(\beta',\beta]\in\pot O_s\setminus\pot O_r$ then potentially
  $\alpha\in(\beta',\beta]$. In that case $\min(N\setminus\alpha)\in
  \pot S_r\cap\ms N\subset\pot S_s$ by (6b) in $r$. But 
  $\min(N\setminus\alpha)\in(\beta',\beta]$ which is in a direct
  contradiction with (7) in $s$.

  {\it Clause (8):} start with some $(\beta',\beta]\in\pot O_r\setminus
  \pot O_s$ and $M\in\pot M_s\setminus\pot M_r$. The interval is
  disjoint from $N$, hence it is disjoint from $M\subset N$.
 
  Now consider arbitrary $M\in\pot M_r\setminus\pot M_s$ and
  $(\beta',\beta]\in\pot O_s\setminus\pot O_r$. Suppose that $(\beta',
  \beta]\cap M\ne\emptyset$ but $\{\beta',\beta\}\not\subset M$. If
  $\beta\ge\sup(M\cap N)$ then, by (9) in $r$, there is an $x$ from
  the $N$-fence for $M$ in the interval. But $x\in\pot S_r\cap\ms N
  \subset\pot S_s$, a contradiction. However, if $\beta<\sup(M\cap N)$
  then, because $\{\beta',\beta\}\not\subset M$, $M\cap N$ is not an
  initial segment of $N$. This means that $M\cap N\in\ms N$. But then
  $(\beta',\beta]$ is inconsistent with $M\cap N$ in $s$, again a
  contradiction.

  {\it Clause (9):} consider arbitrary models $M\in\pot M_r\setminus
  \pot M_s$ and $M'\in\pot M_s\setminus\pot M_r$. Let $\delta:=\sup
  (N\cap M)$, $\delta':=\sup(M'\cap M)=\sup(M'\cap N\cap M)\le\delta$ 
  and $\varepsilon:=\min(N\setminus\delta)$. Let us first establish
  the proper connection between their intersection and each of them.

  First assume that $M\cap N\in\ms N$. Then $M\cap M'=M\cap N\cap M'$
  is in the same correspondence with both $\ms M$ and $\ms M'$ as
  $(M\cap N)\cap M'$ is with $\ms M\cap\ms N$ and $\ms M'$
  respectively. Notice, that $M\cap N$ is an initial segment of $M$.
  On the other hand, if $M\cap N$ is an initial segment of $N$ then
  $M\cap M'$ is also an initial segment of $M'\subset N$. At the same
  time we can prove that $M\cap M'\in\ms M$, as long as we can prove
  that $M\cap M'\in[\delta]^\omega\cap\ms N$, because then we can envoke 
  the fact that $[\delta]^\omega\cap\ms N\subset\ms M$. It is clear that
  $M\cap M'\in[\delta]^\omega$, since $\delta'\le\delta$. If $M\cap M'=
  M'\cap\varepsilon\in\ms N$, then we are done. To that end, first let 
  $\lambda\in M\cap M'$. Then $\lambda\in M\cap N$, hence $\lambda<
  \delta\le\varepsilon$. Therefore $M\cap M'\subset M'\cap\varepsilon$. 
  For the other direction assume that $\lambda\in M'\cap\varepsilon$.
  Then $\lambda\in N\cap\varepsilon$, hence $\lambda\le\delta$. If
  $\lambda=\delta$ then $\delta\in N$, which cannot happen because
  $M\cap N\not\in\ms N$. Therefore $\lambda\in N\cap\delta\subset M$.

  Let us now establish the existence of fences. The $M$-fence for $M'$
  is some subset of the union of the $M$-fence for $N$ and (if 
  $M\cap N\in\ms N$) the $M\cap N$-fence for $M'$, because 
  $\{\min(M\setminus\lambda)\mid\lambda\in M',\,\lambda>
  \sup(M\cap M')\}=\{\min(M\setminus\lambda)\mid\lambda\in N,\,
  \lambda>\sup(M\cap N)\}\cup\{\min((M\cap N)\setminus\lambda)\mid
  \lambda\in M',\,\sup(M\cap M')<\lambda\le\sup(M\cap N)\}$, and this
  union is a subset of $\pot S_r\cup\pot S_s=\pot S_t$.
  
  Let $x\subset\pot S_r\cap\ms N$ be the $N$-fence for $M$ and $y\subset
  \pot S_s$ the $M'$-fence for $M\cap N$ (if $M\cap N\in\ms N$). Then the
  $M'$-fence for $M$ is some subset of $\{\min(M'\setminus\lambda)\mid
  \lambda\in x\cup y\cup\{\min(N\setminus\sup(M\cap N))\}\}\subset\pot
  S_s$.
\theend
\end{prooff}

\section{Preservation of $\omega_2$}

We have thus far proved that forcing with $P$ preserves $\omega_1$. We
also need $\omega_2$ to be preserved. For that purpose we use a weaker
version of closedness, which was also used in~\cite{Mit09}.

\begin{df}
  Assume that forcing notion $P$ preserves cardinals $<\delta$. $P$ is
  {\em $\delta$-presaturated} if for every $A\subset V$, $A\in V[G]$,
  with $|A|^{V[G]}<\delta$, there exists $A'\in V$ such that
  $|A'|^V<\delta$ and $A'\supset A$.
\end{df}

Notice that since $P$ preserves cardinals below $\delta$, $|A'|^V=
|A'|^{V[G]}$ as soon as $|A'|^V<\delta$. Hence we can omit the superscript.

\begin{prop}
  Suppose $\delta$ is a regular cardinal in $V$. If $P$ is
  $\delta$-presaturated then $P$ preserves $\delta$.
\end{prop}

\begin{prooff}
  Suppose for contradiction that $f:\mu\to\delta$ is a bijection in
  $V[G]$ for some $\mu<\delta$. For $\alpha<\mu$, let $A'_\alpha\in V$
  be such that $|A'_\alpha|<\delta$ and $f[\alpha+1]\subset A'_\alpha$.
  Define a function $g:\mu\to\delta$ by $g(\alpha):=\sup(A'_\alpha\cap
  Ord)<\delta$ for $\alpha<\mu$. Then $g\in V$ and $f(\alpha)\le g(\alpha)$
  for all $\alpha<\mu$. Hence $g$ is cofinal, and we get a contradiction.
\theend
\end{prooff}

\begin{lemma}
  Let $\delta$ be a cardinal regular in $V$ such that $P$ preserves
  cardinals below $\delta$. Suppose that for every collection $\pot A$
  of fewer than $\delta$ antichains there exists a dense set
  $\ms D\subset P$ such that for  every $p\in\ms D$, the set 
  $\{q\in\bigcup\pot A\mid p$ and $q$ are compatible$\}$ has size less
  than $\delta$. Then $P$ is $\delta$-presaturated.
\end{lemma}

\begin{prooff}
  Suppose $A\subset V$ and $|A|^{V[G]}<\delta$. Let $p\in G$ be a
  condition such that $p\forces{|\undertilde{A}|<\delta}$. Therefore
  $p\forces{$there exists $\mu<\delta$ such that $|\undertilde{A}|=
  \mu}$. Let $p_0\ge p$, $g$ and $\mu^*<\delta$ be such that 
  $p_0\forces{\undertilde{g}:\mu^*\to\undertilde{A}$ is a
  bijection$}$. For each $\alpha<\mu^*$ let $\pot A_\alpha$ be a
  maximal antichain of conditions in the set $\{q\mid(q\ge p_0\land q$ 
  decides $g(\alpha))\lor q\bot p_0\}$. Hence $\pot A_\alpha$ is a
  maximal antichain above $p_0$.
  
  Define $\pot A:=\{\pot A_\alpha\mid\alpha<\mu^*\}$. Let $\ms D$ be
  a dense set guaranteed by the assumption, and let $p_1\in\ms D$,
  $p_1\ge p_0$. Then the set $X:=\{q\in\bigcup_{\alpha<\mu^*}
  \pot A_\alpha\mid q$ is compatible with $p_1\}$ has size $<\delta$. 
  Let $\Gamma:=\{\beta\mid$ there exist $q\in X$ and $\alpha<\mu^*$
  such that $q\forces{\undertilde{g}(\alpha)=\beta}\}$. 
  Consider an arbitrary $\alpha<\mu^*$. Since $\pot A_\alpha$ is a
  maximal antichain there exists some $q\in\pot A_\alpha$, compatible
  with $p_1$, such that $q$ decides $g(\alpha)$. Hence there exists
  $\beta$ such that $q\forces{\undertilde{g}(\alpha)=\beta}$, and
  therefore $\beta\in\Gamma$. Let $r$ be a common upper bound for $q$
  and $p_1$. Then $r\forces{\undertilde{g}(\alpha)=\beta}$, and since
  $r\ge p_0$, $p_0\forces{$there exists $\beta\in\Gamma$ such that
  $\undertilde{g}(\alpha)=\beta}$. It follows that $p_0\forces{
  \undertilde{g}(\alpha)\in\Gamma}$, so $p_0\forces{\undertilde{g}
  [\mu^*]=\undertilde{A}\subset\Gamma}$. Since $\delta$ is regular in
  $V$, we have $|\Gamma|<\delta$, therefore $p\forces{$there exists
  $A'\in V$, $\undertilde{A}\subset A'$ and $|A'|<\delta}$.
\theend
\end{prooff}

The next lemma shows that $\delta$-presaturation is, in fact, a
generalization of properness to cardinals above $\omega_1$.

\begin{lemma}
  Let $\delta$ be a cardinal regular in $V$ such that $P$ preserves
  cardinals below $\delta$. Suppose that $\theta$ is a large enough
  cardinal, and that for stationarily many models $\ms N$ of size
  $<\delta$ in $[H_\theta]^{<\delta}$ with $P\in\ms N$, and for each
  $p\in P\cap\ms N$, there exists an $\ms N$-generic extension
  $q\ge p$. Then $P$ is $\delta$-presaturated.
\end{lemma}

\begin{prooff}
  Suppose $A\subset V$ and $\mu:=|A|^{V[G]}<\delta$. Let $f$ be a
  function and let $p\in G$ be a condition such that
  $p\forces{\undertilde{f}:\mu\to\undertilde{A}\text{ is onto}}$. 
  Define $\mf N:=\{\ms N\prec H_\theta\mid|\ms N|<\delta,\,\{f,A,p,P\}
  \cup\mu\subset\ms N\}$. Consider some $\ms N\in\mf N$. Let $q\ge p$
  be a generic extension. Then for every $\xi<\mu$, the set $\ms
  D_\xi:=\{r\in\ms N\mid r$ decides $f(\xi)\}\in\ms N$ is dense above $q$.
  Hence $q\forces{\ms D\cap G\cap\ms N\ne\emptyset}$. Therefore $q$
  forces that there exist $r_\xi\in G\cap\ms N$ and $x_\xi\in\ms N$
  such that $r_\xi\forces{\undertilde{f}(\xi)=x_\xi}$. It follows that
  $q\forces{\undertilde{A}\subset\ms N}$, so $p\forces{$there exists
  $A'\in V$, $\undertilde{A}\subset A'$ and $|A'|<\delta}$, $A'$ being
  the model $\ms N$.
\theend
\end{prooff}

\begin{prop}\label{PropPresaturation}
  $P$ is $\omega_2$-presaturated.
\end{prop}

Since presaturation is a generalization of properness, the proof will
be very similar to the proof of properness. Actually, it will be
slightly easier, because we will not work with arbitrary models of size
$\omega_1$ but only with such models that are in a way transitive below
$\omega_2$.

\begin{df}
  Let $\theta>\omega_2$ be some cardinal. Define $\mf M_2:=
  \{\ms M\prec H_\theta\mid |\ms M|=\omega_1$, $\omega_1\cup\{\ms E\}
  \subset\ms M$, $[\ms M]^\omega\subset\ms M\}$.
\end{df} 

Recall that we have assumed CH so that the set $\mf M_2$ is stationary
in $[H_\theta]^{<\omega_2}$. If $\ms M\in\mf M_2$ then $\ms M\cap\omega_2$
is some ordinal $\delta_M\in\omega_2$, since $\omega_1\subset\ms M$
(see~\cite{JusWee97}). Additionally, if $A\in\ms M$ and
$|A|\le\omega_1$ then $A\subset\ms M$.

To prove Proposition~\ref{PropPresaturation}, we first isolate a lemma
which is an analogue of Lemma~\ref{LemmaR}.

\begin{lemma}
  Let $\ms N\in\mf M_2$, and $r\in P$ such that $\delta_N\in\pot S_r$.
  Then $r_{\ms N}$ defined by $\pot F_{r_{\ms N}}:=\pot F_r\cap\ms N$,
  $\pot S_{r_{\ms N}}:=(\pot S_r\cap\ms N)\cup\{\sup(M\cap N)\mid M\in
  \pot M_r\setminus\ms N\}$, $\pot O_{r_{\ms N}}:=\pot O_r\cap\ms N$,
  $\pot M_{r_{\ms N}}:=\{M\cap N\mid M\in\pot M_r\}$, is a condition
  in $P\cap\ms N$.
\end{lemma}

\begin{prooff}
  Notice that, since $M\cap N$ is a countable subset of $\ms N$, it is
  in $\ms N$. Hence, $r_{\ms N}\in\ms N$. Also, if $M\in\ms N$ then
  $M=M\cap N\in\pot M_{r_{\ms N}}$. It is also of some importance that
  by clause (6a), $\pot D_{r_\ms N}=\pot D_r\cap\ms N$.

  Now we have to show, just as in the proof of Lemma~\ref{LemmaR},
  that $r_{\ms N}$ is a condition. By Lemmas~\ref{LemmaIntersection1}
  and~\ref{LemmaIntersection2}, $\ms M\cap\ms N\prec H_{\omega_2}$. At
  the same time, $\ms M\cap\ms N\in\mf M_0$,
  therefore $M\cap N$ can be justifiably added to $\pot M_{r_{\ms N}}$. 
  Since $\delta_N\in\pot S_r$, it follows that $\sup(M\cap N)\in
  \pot D_r\cap\ms N=\pot D_{r_{\ms N}}$ by clause (6b). We can safely
  put it into $\pot S_{r_{\ms N}}$ without violating clause (5),
  because we can apply Lemma~\ref{LemmaLimPoint2}.

  Just as in the proof of Lemma~\ref{LemmaR}, we only have to pay
  attention to clauses (6) and (9) in conjecture with a model of the
  form $M\cap N$. 

  If $\alpha\in M\cap N$ then $C_\alpha\in\ms M$ by (6a) in $r$. If
  $\cf(\alpha)=\omega$ then $C_\alpha\in\ms N$, because $\ms N$ is
  countably closed. If $\cf(\alpha)=\omega_1$ then $C_\alpha\in\ms N$,
  because $C_\alpha=C\setminus\beta$ for some $C\in\ms E_\alpha\subset
  \ms N$ and $\beta\in\pot D_r\cap\alpha\subset\ms N$. In any case,
  $C_\alpha\in\ms M\cap\ms N$. The danger with clauses (6c) and (6d)
  is if $\alpha\in\ms N$ is in a gap of $M\cap N$ and $\alpha\in
  \Lim(C_\gamma)$ for some $\gamma\in\pot D_r\setminus\ms N$. But then
  $\gamma>\delta_N\in\pot S_r$, hence $\alpha\not\in\Lim(C_\gamma)$. 

  Let $M_1':=M_1\cap N$ and $M_2':=M_2\cap N$ be two models from 
  $\pot M_{r_{\ms N}}$. Suppose that $M_1\cap M_2$ is an initial segment
  of $M_1$ and $A\in[\lambda']^\omega\cap\ms M_1'$ where, $\lambda':=
  \sup(M_1'\cap M_2')\le\sup(M_1\cap M_2)=:\lambda$. Then $A\in
  [\lambda]^\omega\cap\ms M_1\cap\ms N\subset\ms M_2\cap\ms N=\ms M_2'$.

  Suppose now that $M_1\cap M_2\in\ms M_1$. Then $\lambda'\in\ms M_1$,
  hence $M_1'\cap M_2'=M_1\cap M_2\cap\lambda'\in\ms M_1$. On the other
  hand, $M_1'\cap M_2'$ is a countable subset of $\ms N$, hence it is
  in $\ms N$. Therefore, $M_1'\cap M_2'\in\ms M_1\cap\ms N=\ms M_1'$.

  The $M_1'$-fence for $M_2'$ is the $M_1$-fence for $M_2$ intersected
  with $N$. The $M_2'$-fence for $M_1'$ is obtained in the same way.
\theend
\end{prooff}

\begin{prooff}[of Proposition~\ref{PropPresaturation}]
  Let $\theta$ be a large enough cardinal. Pick $\ms N\in\mf M_2$ and
  $p\in P\cap\ms N$. We extend $p$ to $q$ by putting $\delta_N$ into
  both $\pot D_p$ and $\pot S_p$. For the corresponding club
  $C_{\delta_N}$ we take $C\setminus\max(\pot D_p)$ for some $C\in
  \ms E_{\delta_N}$.\marginpar{\tiny Possible complication: is $C_\delta$
  in $\ms E_\delta$?}
  Clearly, $q\in P$. We will prove that $q$ is an 
  $\ms N$-generic extension of $p$.

  Suppose $r$ is an arbitrary extension of $q$, and let $r_{\ms N}$ be
  as given by the previous lemma. For a fixed dense set $\ms D\subset P$, 
  $\ms D\in\ms N$, extend $r_{\ms N}$ to $s\in\ms D$. Then $s\in\ms N$. 
  As with properness, we will prove clause by clause of
  Definition~\ref{DefFor} that $t:=(\pot F_r\cup\pot F_s,\,\pot S_r\cup
  \pot S_s,\,\pot O_r\cup\pot O_s,\,\pot M_r\cup\pot M_s)$ is a condition.
  Fortunately, less effort will have to be invested, because $N\cap 
  H_{\omega_2}$ has no gaps. 

  {\it Clauses (1), (2) and (3)} need no comments.

  {\it Clause (4):} suppose that $\alpha\in\pot D_r\setminus\pot D_s$
  and $\beta\in\pot D_s\setminus\pot D_r$. Then $C_\alpha\cap\ms N$ is
  a finite set because $\delta_N\in\pot S_r$, and $C_\beta\subset\ms N$,
  hence $\Lim(C_\alpha)\cap\Lim(C_\beta)=\emptyset$. 

  {\it Clause (5):} if $\alpha\in\pot D_r\setminus\pot D_s$ and
  $\sigma\in\pot S_s\setminus\pot S_r$ then $C_\alpha\cap\sigma\subset
  C_\alpha\cap\delta_N$, which is a finite set.

  {\it Clause (6):}  if $\alpha\in\pot D_r\setminus\pot D_s$ and
  $M\in\pot M_s\setminus\pot M_r$ then there is nothing to
  prove. Suppose now that $\alpha\in\pot D_s\setminus\pot D_r$ and
  $M\in\pot M_r\setminus\pot M_s$. Then $M\cap N\in\pot M_s$ and the
  compatibility between $\alpha$ and $M\cap N$ in $s$ is transfered to
  the compatibility between $\alpha$ and $M$. The only potential
  problem would be, just as in the proof of the previous lemma, if 
  $\alpha\ge\sup(M\cap N)$ and $\alpha\in\Lim(C_\gamma)$ for some
  $\gamma\in\pot D_r\setminus\pot D_s$. But as we saw, $\Lim(C_\gamma)
  \cap\ms N=\emptyset$, because $C_\gamma\cap\delta_N$ is a finite set. 

  {\it Clause (7):} if $\alpha\in\pot D_r\setminus\pot D_s$ and
  $(\beta',\beta]\in\pot O_s\setminus\pot O_r$ then $(\beta',\beta]
  \subset\ms N$, hence $\alpha\not\in(\beta',\beta]$. Suppose now that
  $\alpha\in\pot D_s\setminus\pot D_r$ and $(\beta',\beta]\in\pot O_r
  \setminus\pot O_s$. Since $\delta_N\in\pot S_r$, we have $(\beta',\beta]
  \cap\ms N=\emptyset$, hence $\alpha\not\in(\beta',\beta]$.

  {\it Clause (8):} if $M\in\pot M_s\setminus\pot M_r$ and 
  $(\beta',\beta]\in\pot O_r\setminus\pot O_s$ then $(\beta',\beta]\cap 
  \ms M=\emptyset$. Consider now some $M\in\pot M_r\setminus\pot M_s$
  and $(\beta',\beta]\in\pot O_s\setminus\pot O_r$. Then $(\beta',\beta]$ 
  and $M\cap N$ are compatible in $s$. If $(\beta',\beta]\in\ms M\cap
  \ms N$ then $(\beta',\beta]\in\ms M$. If $(\beta',\beta]\cap(\ms M
  \cap\ms N)=\emptyset$ then $(\beta',\beta]\cap\ms M=
  ((\beta',\beta]\cap\ms N))\cap\ms M=\emptyset$. 

  {\it Clause (9):} consider two models $M\in\pot M_r\setminus\pot M_s$ 
  and  $M'\in\pot M_s\setminus\pot M_r$. Then $M\cap N$ and $M'$ are
  compatible in $s$. If $M\cap N\cap M'\in\ms M\cap\ms N$ then $M\cap M'= 
  M\cap M'\cap N\in\ms M$. Now suppose that $M\cap N\cap M'$ is an
  initial segment of $M\cap N$ and let $A\in[\lambda]^\omega\cap
  \ms M$, where $\lambda:=\sup(M\cap M')=\sup(M\cap N\cap M')\in N$. 
  Then $A\in\ms N$, since $\ms N$ is countably closed. Hence 
  $A\in[\lambda]^\omega\cap\ms M\cap\ms N\subset\ms M'$.

  On the other hand, if $M\cap N\cap M'\in\ms M'$ then $M\cap M'= 
  M\cap M'\cap N\in\ms M'$. However, if $M\cap N\cap M'$ is an
  initial segment of $M'$ and $A\in[\lambda]^\omega\cap\ms M'\subset
  \ms M\cap\ms N$, then $A\in\ms M$, hence $M\cap M'$ is an initial
  segment of $M'$.
\theend
\end{prooff}

\begin{posl}
  Forcing $P$ preserves cardinals.
\end{posl}

\begin{prooff}
  $P$ has the $\cc{\omega_3}$ because, assuming $2^{\omega_1}=\omega_2$,
  $|P|=\omega_2$. Hence it preserves cardinals $\ge\omega_3$. It 
  preserves $\omega_2$ because it is $\omega_2$-presaturated. And it
  preserves $\omega_1$ because it is proper.
\theend
\end{prooff}

\begin{df}
  Let $G\subset P$ be a generic set. Define $\pot F:=\bigcup_{p\in G}
  \pot F_p$, and $\pot C:=\dom(\pot F)$.
\end{df}

We have no reason to believe that $\pot C=\Lim(\omega_2)$. The usual
density argument fails in this case. Namely, we cannot extend a given
$p\in P$ with $\alpha\not\in\pot D_p$ to $q\in P$ such that $\alpha\in
\pot D_q$. At least not for any $p\in P$ and $\alpha<\omega_2$. Suppose
$p=(\{(\omega_1+\omega,C)\},\emptyset,\emptyset,\{M\})$ and
$C=\omega+1\cup[\omega_1,\omega_1+\omega)\in\ms M$. If we want to
extend $\pot D_p$ to $\alpha\in[\delta_M,\omega_1]$ then we have to add
$\omega_1$ to $\pot S_p$ for clause (6b) to hold. But then $\omega_1$
and $C$ violate clause (5).

\begin{prop}
  $\pot C$ is unbounded in $\omega_2$.
\end{prop}

\begin{prooff}
  Define $\ms D_\alpha:=\{p\in P\mid\max(\pot D_p)>\alpha\}$ for
  $\alpha<\omega_2$. Consider arbitrary $p\in P$ and assume that
  $p\not\in\ms D_\alpha$. Now let $\alpha':=\sup(\pot D_p\cup\bigcup
  \pot O_p\cup\bigcup\pot M_p)<\omega_2$ and $q:=(\pot F_p\cup\{\alpha'+
  \omega,(\alpha',\alpha'+\omega)\},\,\pot S_p,\,\pot O_p,\,\pot M_p)$. 
  Clearly, $q\in P$, $q\ge p$ and $q\in\ms D_\alpha$, hence
  $\ms D_\alpha$ is dense in $P$ for every $\alpha<\omega_2$. It
  follows that $\pot C$ is unbounded in $\omega_2$. 
\theend
\end{prooff}

To prove that $\pot C$ is closed, we need the following lemma.

\begin{lemma}
  Suppose that $p\in P$ is such that $\alpha\not\in\pot D_p$. There
  exists an extension $q\ge p$ such that there is some interval
  $(\beta',\beta]\in\pot O_q$ with $\alpha\in(\beta',\beta]$. 
\end{lemma}

\begin{prooff}
  Assume that $\alpha\not\in(\beta',\beta]$ for every
  $(\beta',\beta]\in\pot O_p$. Define $A:=\{M\in\pot M_p\mid
  \alpha\not\in M,\ \sup(M\cap\alpha)=\alpha\}$ and 
  $B:=\{M\in\pot M_p\mid\alpha\in M\}$. Notice that by clause (3) of
  Definition~\ref{DefFor}, $\alpha<\sup(M)$ for every $M\in A$.

  First assume that $B=\emptyset$. If $A=\emptyset$ then define
  $\beta':=\sup((\pot D_p\cup\bigcup\pot M_p)\cap\alpha)<\alpha$ and
  $\beta:=\alpha$. Obviously $(\beta',\beta]\cap\pot D_p=\emptyset$
  and $(\beta',\beta]\cap M=\emptyset$ for every $M\in\pot M_p$. 

  If $A\ne\emptyset$ then let $\beta'$ be any ordinal in $\bigcap A$
  strictly between $\sup((\pot D_p\cup\bigcup(\pot M_p\setminus A))
  \cap\alpha)$ and $\alpha$. Such an ordinal exists because $\alpha$
  is a supremum of some block in every $M\in A$. Namely, if $M,M'\in A$ 
  have interchanging blocks with supremum $\alpha$ then there has to be
  infinitely many such blocks, but that cannot happen neither below
  $\sup(M\cap M')$ nor above it (see Remark~\ref{CompatibleModels}). 
  Therefore only one model from $A$ can have $\alpha$ as a supremum of
  some sequence of blocks rather than a supremum of a certain block. 
  To be precise, there can be more such models but the sequence must
  be the same from some point on. In this case we can still find
  $\beta'$ as described.

  On the other hand, if $M,M'\in A$ then $\min(M\setminus\alpha)=
  \min(M'\setminus\alpha)$. That is true because if $\alpha=\sup(M\cap
  M')$ then $\alpha\in\pot D_p$. But if $\alpha<\sup(M\cap M')$ and
  if, for instance, $\alpha':=\min(M\setminus\alpha)<\min(M'\setminus
  \alpha)$ then $M\cap M'$ is not an initial segment of $M$, hence it
  must be an element of $M$, which means that $\alpha=\sup((M\cap M'))
  \cap\alpha')\in M$. However, that contradicts the original
  assumption that $\alpha\not\in M$. Define $\beta:=\min(M\setminus
  \alpha)$ for some (i.e.~every) $M \in A$. Interval $(\beta',\beta]$ is
  an element of every $M\in A$. If $\gamma\in(\beta',\beta]$ for some
  $\gamma\in\pot D_p$ then $\alpha\in\pot D_p$ by clause (6b) of
  Definition~\ref{DefFor}. If $(\beta',\beta]\cap M'\ne\emptyset$ for
  some $M'\in\pot M_p\setminus A$ then we have to consider two
  distinct possibilities. If it happens above the $\sup(M\cap M')$ for
  some $M\in A$ then a block of $M'$ can only lie inside a gap of $M$,
  which means that $\beta$ is in the $M$-fence for $M'$, and, again by
  (6b), $\alpha\in\pot D_p$. 

  If it occurs below the $\sup(M\cap M')$ for some $M\in A$ then it
  can only happen if $(\beta',\beta]\cap M'=\{\beta\}$. In this case
  $(\beta',\beta]$ and $M'$ are not compatible. Let $A':=\{M'\in\pot
  M_p\setminus A\mid\beta\in M'\}$. We have to extend the interval
  $(\beta',\beta]$ downwards so that $\beta'\in M'$ for every $M'\in A'$ 
  but $\beta'>\sup((\pot D_p\cup\bigcup(\pot M_p\setminus(A\cup A')))
  \cap\alpha)$. If we can find such $\beta'$ then this new interval
  will still be compatible with every $M\in\pot M_p\setminus A'$, and
  it will obviously be compatible with every $M'\in A'$. We can be sure
  that there is no $\gamma\in\pot D_p$ in the gap between $\sup(M'\cap
  \beta)$ and $\beta$,  because otherwise $\beta\in\pot D_p$, hence
  $\alpha\in\pot D_p$. So this takes care of $\pot D_p$. Can we find
  $\beta'$ so that the new interval does not intersect any $M''\in\pot
  M_p\setminus(A\cup A')$? Suppose $M''$ is the reason that the answer
  is negative. Then $\sup(M''\cap\beta)\ge\sup(M'\cap\beta)$. If
  $\sup(M''\cap\beta)\ge\sup(M'\cap M'')$ then $\beta$ is in the
  $M'$-fence for $M''$ and we derive a contradiction. If $\sup(M''\cap 
  \beta)<\sup(M'\cap M'')$ then either we get a contradiction by
  Remark~\ref{CompatibleModels} if $\sup(M''\cap\beta)>\sup(M'\cap\beta)$,
  or, if $\sup(M''\cap\beta)=\sup(M'\cap\beta)$, we can conclude that
  $M'\cap M''\in M'$ (because $\beta\in M'\setminus M''$), hence
  $\sup(M'\cap\beta)=\sup((M'\cap M'')\cap\beta)\in M'$, again a
  contradiction. We should check that we can find $\beta'$ so that it
  lies in every $M\in A$. If $M\in A$ is a counterexample then the
  whole block of $M$ lies inside of a gap in some $M'\in A'$. It
  happens below the $\sup(M\cap M')$, as $\beta\in M\cap M'$, hence we
  can use Remark~\ref{CompatibleModels} to get a contradiction.

  Now assume that $B\ne\emptyset$ and $A=\emptyset$. If
  $\sup(M\cap\alpha)=\alpha$ for every $M\in B$ then, similarly as
  above, let $\beta'$ be any ordinal in $\bigcap B$ strictly between
  $\sup((\pot D_p\cup\bigcup(\pot M_p\setminus B))\cap\alpha)$ and
  $\alpha$, and $\beta:=\alpha$. Then $(\beta',\beta]\in M$ for every
  $M\in B$ while it is disjoint from $\pot D_p$ and every $M'\in\pot
  M_p\setminus B$. If there is some $M\in B$ such that $\sup(M\cap
  \alpha)<\alpha$ then we construct $(\beta',\beta]$ in the very same
  way, but we have to check that $\sup((\pot D_p\cup\bigcup(\pot M_p
  \setminus B))\cap\alpha)<\sup(M\cap\alpha)$ so that $\beta'$ can be
  taken from $M$. If there is some $\gamma\in\pot D_p$ between $\sup(
  M\cap\alpha)$ and $\alpha$ then $\alpha\in\pot D_p$. The same is
  true if there is a block of some $M'\in\pot M_p\setminus B$ inside
  the gap of $M$ just below $\alpha$ and $\sup(M\cap\alpha)\ge\sup(M
  \cap M')$. However, if $\sup(M\cap\alpha)<\sup(M\cap M')$ then there
  can be nothing from $M'$ inside the gap of $M$ below $\alpha$,
  because, as several times before, that would make it impossible for
  $M\cap M'$ to be an initial segment of either $M$ or $M'$. This
  follows from the fact that $\alpha\in M\setminus M'$.

  Finally assume that $B\ne\emptyset$ and $A\ne\emptyset$. If $\sup(M
  \cap\alpha)=\alpha$ for every $M\in B$ then let $\beta'$ be any ordinal
  in $\bigcap(B\cup A)$ strictly between $\sup((\pot D_p\cup\bigcup(\pot
  M_p\setminus(B\cup A)))\cap\alpha)$ and $\alpha$. To define $\beta$,
  first recall that if $M,M'\in A$ then $\min(M\setminus\alpha)=\min
  (M'\setminus\alpha)$. Additionally, if $M\in B$ and $M'\in A$ then
  $\min(M'\setminus\alpha)\in M$. This is true because $\alpha<\sup(M
  \cap M')$, hence if it were not true then $M\cap M'$ would not be an
  initial segment of either $M$ or $M'$. Notice here that if $\alpha= 
  \sup(M\cap M')$ then $\alpha\in\pot D_p$. Define $\beta:=\min(M'
  \setminus\alpha)$ for some $M'\in A$. Then $(\beta',\beta]\in M$ for
  every $M\in B\cup A$. It is also disjoint from $\pot D_p$ and every
  $M'\in\pot M_p\setminus(B\cup A)$, because otherwise $\alpha\in\pot
  D_p$. The only possible exception is if $(\beta',\beta]\cap M'=
  \{\beta\}$ for some $M'\in\pot M_p\setminus(B\cup A)$. This
  situation is dealt with in the same manner as in the case if
  $B=\emptyset$ and $A\ne\emptyset$. If $\sup(M\cap\alpha)<\alpha$ for
  some $M\in B$ while at the same time $\sup(M'\cap\alpha)=\alpha$ for
  some $M'\in A$ then it has cofinality $\omega_1$ (see
  Lemma~\ref{LemmaTopGap}) and $\omega$. 

  In all the cases we were able to find an interval $(\beta',\beta]$
  compatible with everything in $p$, such that $\alpha\in(\beta',\beta]$.
  Hence, $q:=(\pot F_p,\pot S_p,\pot O_p\cup\{(\beta',\beta]\},\pot M_p)$
  is the desired extension of $p$.
\theend
\end{prooff}

\begin{prop}
  $\pot C$ is closed in $\omega_2$.
\end{prop}

\begin{prooff}
  Suppose for contradiction that $p\in G$ is such that
  $p\forces{\alpha\in\Lim(\pot C)$ but $\alpha\not\in\pot C}$ for some
  $\alpha<\omega_2$. Then $\alpha\not\in\pot D_p$. Let $q$ be the
  extension given by previous lemma. But then
  $q\forces{\alpha\not\in\Lim(\pot C)}$, which contradicts the fact
  that $p\forces{\alpha\in\Lim(\pot C)}$.
\theend
\end{prooff}

What we have created might not be a $\Box_{\omega_1}$ sequence, but
the next lemma shows that we can now extend our square-like sequence to
the whole $\Lim(\omega_2)$.

\begin{lemma}\label{LemmaSquareExt}
  Let $\kappa$ be a regular cardinal $>\omega$. Suppose that $\pot C
  \subset\Lim(\kappa^+)$ is a club of $\kappa^+$ and $\langle C_\alpha
  \mid\alpha\in\pot C\rangle$ is a square-like sequence. Then there
  exists a square sequence on $\kappa^+$.
\end{lemma}

\begin{prooff}
  The idea is to throw away every ordinal which is not in $\pot C$,
  effectively making $\pot C$ equal to $\kappa^+$. In fact, keeping
  only limit points of $\pot C$ will suffice. Thus, let $\pot E:=
  \Lim(\pot C)\setminus\{\kappa^+\}$. $\pot E$ is stil a club of
  $\kappa^+$. For every $\alpha\in\Lim(\pot E)\setminus\{\kappa^+\}$
  define $D_\alpha:=C_\alpha\cap\pot E$. Since $\pot E\cap\alpha$ is a
  club in $\alpha$ for every $\alpha\in\Lim(\pot E)\setminus\{\kappa^+\}$,
  $D_\alpha$ is a club in $\alpha$. Suppose that $\beta\in
  \Lim(D_\alpha)$ for some $\beta<\alpha$. Then $\beta$ is a limit
  point of both $\pot E$ and $C_\alpha$, and $D_\beta=C_\beta\cap\pot E=
  C_\alpha\cap\beta\cap\pot E=D_\alpha\cap\beta$. Also, if $\cf(\alpha)
  <\kappa$ then $|D_\alpha|<\kappa$. Hence, $\langle D_\alpha\mid\alpha
  \in\Lim(\pot E)\setminus\{\kappa^+\}\rangle$ is a square-like sequence.

  Let $\{\gamma_i\mid i<\kappa^+\}$ be an increasing enumeration of
  $\pot E$. For $i\in\Lim(\kappa^+)$ define $E_i:=\{j<i\mid\gamma_j\in
  D_{\gamma_i}\}=\gamma^{-1}[D_{\gamma_i}]$. It is a club in $i$
  because $\gamma$ is a continuous function. Let us prove that
  $\langle E_i\mid i\in\Lim(\kappa^+)\rangle$ is a square sequence. If
  $i<j$ and $i\in\Lim(E_j)$ then $\gamma_i\in\Lim(D_{\gamma_j})$.
  Hence, $D_{\gamma_i}=D_{\gamma_j}\cap\gamma_i$. Therefore, $E_i=
  \gamma^{-1}[D_{\gamma_i}]=\gamma^{-1}[D_{\gamma_j}\cap\gamma_i]=
  \gamma^{-1}[D_{\gamma_j}]\cap i=E_j\cap i$. If $\cf(i)<\kappa$ then
  $\cf(\gamma_i)<\kappa$, hence $|E_i|=|D_{\gamma_i}|<\kappa$.
\theend
\end{prooff}

\begin{posl}
  $V[G]\models\Box_{\omega_1}$.
\end{posl}

\bibliographystyle{amsplain}
\bibliography{bibliography}

\end{document}